\documentstyle[12pt,fleqn,leqno]{article}
\def\lastrevision{April 29, 1996}
\title{{\LARGE\bf More set-theory\protect\\ {\it around}\protect\\
 the weak Freese-Nation property}\\} 
\author{Saka\'e Fuchino$^{0}$, Lajos Soukup$^{1}$}
\date{\lastrevision}
\newif\ifwrong
\newif\iftesting
\newif\ifcommented 
\let\nc\newcommand
\let\rnc\renewcommand
\let\nth\newtheorem
\let\nenv\newenvironment
\iftesting
\let\Label\label%
\def\label#1{\mbox{}\marginpar{{\tiny #1}}\Label{#1}}%
\fi
\ifcommented
\fi
\nenv{comments}{\ifcommented\begin{footnotesize}}%
{\ifcommented\end{footnotesize}\fi}
\nc{\?}{\iftesting%
\marginpar{$\emptyset\emptyset\emptyset$}%
\fi
\typeout{Something is still incomplete here.}
\errmessage{................................}}
\ifx\kanjiskip\otankonasu
\else\fi
\rnc{\baselinestretch}{1.17}
\nth{Thm}{{\bf Theorem}}
\nth{Cor}[Thm]{{\bf Corollary}}
\nth{Lemma}[Thm]{{\bf Lemma}}
\nth{Prop}[Thm]{{\bf Proposition}}
\nth{Claim}{{\bf Claim}}[Thm]
\nth{Subclaim}{{\bf Subclaim}}[Claim]
\nth{Problem}{{\bf Problem}}
\nth{Ex}[Thm]{{\bf Example}}
\nth{Fact}[Thm]{{\bf Fact}}
\nc{\Thmof}[1]{{Theorem \ref{#1}}}
\nc{\Corof}[1]{{Corollary \ref{#1}}}
\nc{\Propof}[1]{{Proposition \ref{#1}}}
\nc{\Lemmaof}[1]{{Lemma \ref{#1}}}
\nc{\Factof}[1]{{Fact \ref{#1}}}
\nc{\Claimof}[1]{{Claim \ref{#1}}}
\nc{\Exof}[1]{{Example \ref{#1}}}
\nc{\Problemof}[1]{{Problem \ref{#1}}}
\nc{\Thmabove}{{Theorem \number\theThm}}
\nc{\Corabove}{{Corollary \number\theThm}}
\nc{\Propabove}{{Proposition \number\theThm}}
\nc{\Lemmaabove}{{Lemma \number\theThm}}
\nc{\Factabove}{{Fact \number\theThm}}
\nc{\Claimabove}{{Claim \number\theClaim}}
\nc{\Exabove}{{Example \number\theThm}}
\nc{\prf}{{\bf Proof\ \ }}
\nc{\prfof}[1]{{\bf Proof of #1\ \ }}
\nc{\prfofClaim}{\raisebox{-.4ex}{\Large $\vdash$\ \ }}
\newsavebox{\qedbox}\sbox{\qedbox}{
{\unitlength=0.034ex \begin{picture}(40,60)
\put(0,0){\framebox(30,44)[cc]{}}
\put(30,-7){\rule{7\unitlength}{44\unitlength}}
\put(10,-7){\rule{27\unitlength}{7\unitlength}}
\end{picture}}}
\nc{\qed}{\nolinebreak\mbox{}\nolinebreak\hfill%
\ifvmode\mbox{}\hfill\fi\usebox{\qedbox}}
\nc{\qedskip}{\medskip\smallskip}
\nc{\smallqed}%
{\mbox{}\smallskip\hfill\raisebox{-.4ex}{\Large $\dashv$}\\}
\nc{\qedof}[1]%
{\nopagebreak\mbox{}\nolinebreak\hspace*{\fill}{\usebox{\qedbox}{~(#1)}}}
\nc{\Qedof}[1]%
{\nopagebreak\mbox{}\nolinebreak\hspace*{\fill}{\usebox{\qedbox}%
{~(#1~\number\theThm)}}}
\nc{\qedofThm}{\Qedof{Theorem}}
\nc{\qedofCor}{\Qedof{Corollary}}
\nc{\qedofProp}{\Qedof{Proposition}}
\nc{\qedofLemma}{\Qedof{Lemma}}
\nc{\qedofFact}{\Qedof{Fact}}
\nc{\qedofClaim}%
{\nopagebreak\mbox{}\nolinebreak\hfill\raisebox{-.4ex}{\Large $\dashv$ }\nolinebreak%
\mbox{~(Claim~\number\theClaim)}}
\nc{\qedofSubclaim}%
{\mbox{}\hfill\raisebox{-.4ex}{\Large $\dashv$ }\nolinebreak%
\mbox{~(Subclaim~\number\theSubclaim)}}
\nc{\noindentafterqed}{\bigskip\\}
\nc{\indentafterqed}{\bigskip\\\hspace*{\parindent}}
\nc{\xitem}[1]{\item[$(\mbox{\it #1})$]}
\nc{\xitemof}[1]{{$(\mbox{\it #1})$}}
\def\assert#1{\makebox[5ex][l]{\rm #1)}\ignorespaces}
\def\lassert#1{\llap{\makebox[5ex][l]{\rm #1)}}\ignorespaces}
\def\assertof#1{{\rm #1)}}
\nenv{assertion}[1]{\begin{trivlist}
\newbox\assertbox
\dimen255=\textwidth
\setbox\assertbox=\hbox{\hspace*{\parindent}{#1}\hspace{\labelsep}}
\advance\dimen255 by -1\wd\assertbox
\advance\dimen255 by -1ex
\item[]\unhbox\assertbox\hfill
\begin{minipage}[t]{\dimen255}}%
{\end{minipage}\end{trivlist}}
\nenv{subassertion}[1]{\begin{trivlist}
\newbox\assertbox
\dimen255=\textwidth
\setbox\assertbox=\hbox{\hspace*{2\parindent}{#1}\hspace{\labelsep}}
\advance\dimen255 by -1\wd\assertbox
\advance\dimen255 by -1ex
\item[]\unhbox\assertbox\hfill
\begin{minipage}[t]{\dimen255}}%
{\end{minipage}\end{trivlist}}
\nc{\implies}{$\,\Rightarrow\,$}
\nc{\equivto}{$\,\Leftrightarrow\,$}
\nenv{formula}[1]{\\[\abovedisplayskip]${#1}$\hfill$}%
{$\hfill\mbox{}\\[\belowdisplayskip]}
\nc{\restr}%
{\mathop{\hspace{0.01ex}|\hspace*{-0.02ex}{\grave{}}\hspace{0.4ex}}}
\nc{\upper}{\mathop{\hspace{0.1ex}\uparrow\hspace{0.2ex}}}
\nc{\subsetnoneq}%
{\mathrel{\raisebox{-0.8ex}{$\stackrel{\subset}{\scriptstyle\,\not=\,}$}}}
\nc{\Linf}[1]{{\cal L}_{\infty{#1}}}
\nc{\powersetof}[1]{{\cal P}(#1)}
\nc{\psof}[1]{{\cal P}(#1)}
\nc{\cardof}[1]{{\mathopen{\mid}{#1}\mathclose{\mid}}}
\nc{\generatedby}[1]{{\langle#1\rangle}}
\nc{\devides}{^{^{\,}}\mid^{^{\,}}}
\nc{\setof}[2]{\{\,#1\,:\,#2\,\}}
\nc{\smallsetof}[1]{\{\,#1\,\}}
\nc{\orderedsetof}[2]{\langle\,#1\,:\,#2\,\rangle}
\nc{\tuple}[1]{\langle\,#1\,\rangle}
\nc{\continuum}{2^{\aleph_0}}
\nc{\reals}{\mbox{$\rm I_{\!\!}R$}}
\nc{\rationals}{%
\rlap{\hbox{\hspace*{0.35ex}\raise0.35ex\hbox{$\scriptstyle\wr$}}}%
\hbox{\rm Q}}
\nc{\integers}{Z\hspace*{-1.12ex}Z}
\nc{\unitint}{\mbox{$\rm I_{\!\!}I$}}
\nc{\lvalue}{\lbrack\!\lbrack}
\nc{\rvalue}{\rbrack\!\rbrack}
\nc{\Bvalueof}[2]{\lvalue\,#1\,\rvalue^{(#2)}}
\nc{\forces}[2]{\,\|\hspace{-.35ex}\mbox{\sf--}_{\,#1\,}%
\mbox{\rm``}\,#2\,\mbox{\rm''}}
\let\Models\models
\def\models#1{\Models\mbox{\rm``}\,#1\,\mbox{\rm''}}
\nc{\xmbox}[1]{ ${\rm #1}$ }
\nc{\st}{such that}
\nc{\wolog}{without loss of generality}
\nc{\Wolog}{Without loss of generality}
\nc{\wrt}{with respect to}
\nc{\Tfae}{The following are equivalent}
\nc{\tfae}{the following are equivalent}
\nc{\Ba}{Boolean algebra}
\nc{\Bas}{Boolean algebras}
\nc{\cBa}{complete Boolean algebra}
\nc{\cBas}{complete Boolean algebras}
\nc{\po}{partial ordering}
\nc{\pos}{partial orderings}
\nc{\Pos}{Partial orderings}
\nc{\Fr}{{\rm Fr}\,}
\nc{\cof}{\mathop{\rm cof}}
\nc{\cf}{\mathop{\rm cf}}
\nc{\pcf}{\mathop{\rm pcf}}
\nc{\Sub}{{\rm Sub}}
\nc{\pr}{{\rm pr}}
\nc{\proj}{{\rm proj}}
\nc{\wproj}{{\rm wproj}}
\nc{\rc}{{\rm rc}}
\nc{\cm}{{cm}}
\nc{\reg}{{\rm reg}}
\nc{\Card}{{\rm Card}}
\nc{\Ord}{{\rm Ord}}
\nc{\Fn}{{\rm Fn}}
\nc{\dom}{{\rm dom}}
\nc{\rng}{{\rm rng}}
\nc{\Lim}{{\rm Lim}}
\nc{\otp}{{\rm otp}}
\nc{\cohenalg}[1]%
{\mbox{$\,\raisebox{0.05ex}{\small$\wr$}\!\!_{_{\!\!}}\mbox{\rm C}_{#1}$}}
\nc{\cantorsp}[1]{\mbox{}^{#1^{\mbox{}\!}}2}
\nc{\fnsp}[2]{\mbox{}^{#1^{\mbox{}\!}}#2}
\nc{\mapping}[3]{#1:#2\rightarrow #3}
\nc{\pfeil}{\mathrel{{\rightarrow\,}\llap{$\rightarrow$}}}
\nc{\ccc}{{\mbox{\rm ccc}}}
\nc{\ZFC}{{\rm ZFC}}
\nc{\CH}{{\rm CH}}
\nc{\GCH}{{\rm GCH}}
\nc{\MA}{{\rm MA}}
\nc{\MM}{{\rm MM}}
\nc{\MAp}{${\rm MA}^+(\sigma${\it-closed\/}$)$}
\nc{\stick}%
{\mbox{{\hspace{0.4ex}$\mbox{\raisebox{1ex}{$\bullet$}}%
\hspace*{-0.94ex}|$\hspace{0.6ex}}}}
\nc{\Subrc}[1]{\Sub_\rc^{\aleph_0}(#1)}
\nc{\dpr}{\overline{\pr}}
\nc{\calA}{{\cal A}}
\nc{\calB}{{\cal B}}
\nc{\calC}{{\cal C}}
\nc{\calD}{{\cal D}}
\nc{\calF}{{\cal F}}
\nc{\calG}{{\cal G}}
\nc{\calH}{{\cal H}}
\nc{\calM}{{\cal M}}
\nc{\calN}{{\cal N}}
\nc{\calO}{{\cal O}}
\nc{\calP}{{\cal P}}
\nc{\calS}{{\cal S}}
\nc{\calT}{{\cal T}}
\nc{\dota}{{\dot{a}}}
\nc{\dotb}{{\dot{b}}}
\nc{\dotc}{{\dot{c}}}
\nc{\dotf}{{\dot{f}}}
\nc{\dotg}{{\dot{g}}}
\nc{\dotm}{{\dot{m}}}
\nc{\dotr}{{\dot{r}}}
\nc{\dotx}{{\dot{x}}}
\nc{\doty}{{\dot{y}}}
\nc{\dotB}{{\dot{B}}}
\nc{\dotP}{{\dot{P}}}
\nc{\dotQ}{{\dot{Q}}}
\nc{\dotR}{{\dot{R}}}
\nc{\dotS}{{\dot{S}}}
\nc{\dotU}{{\dot{U}}}
\nc{\dotX}{{\dot{X}}}
\nc{\dotY}{{\dot{Y}}}
\nc{\undercirc}[1]{%
\setbox255=\hbox{$#1$}
\dimen255=0.9ex
\advance\dimen255 by \dp255
\rlap{\hbox to\wd255{\hss{\lower\dimen255\hbox{$\scriptstyle\circ$}}\hss}}
\box255
}
\nc{\sundercirc}[1]{%
\setbox255=\hbox{$\scriptstyle#1$}
\dimen255=0.8ex
\advance\dimen255 by \dp255
\rlap{\hbox to\wd255{\hss{\lower\dimen255
\hbox{$\scriptscriptstyle\circ$}}\hss}}
\box255
}
\nc{\multiundercirc}[1]{%
\mathchoice{\undercirc{#1}}{\undercirc{#1}}{\sundercirc{#1}}{\sundercirc{#1}}}
\nc{\acirc}{\multiundercirc{a}}
\nc{\Dag}{^\dagger}
\nc{\deq}{\mathrel{\sqsubseteq}}
\nc{\Deq}{\mathrel{{\sqsubseteq}\hspace{-1.5ex}{\sqsubseteq}}}
\nc{\notDeq}{\mathrel{{\sqsubseteq}\hspace{-1.5ex}{\not\sqsubseteq}}}
\nc{\kappagame}[2]{\calG^{#1}(#2)}%
\nc{\kappalambdagame}[3]{\calG^{#1}_{#2}(#3)}%
\begin{document}
\maketitle
\footnotetext{\ Institut f\"ur Mathematik II, Freie Universit\"at Berlin, 
14195 Berlin, Germany. e-mail:
{\tt fuchino@math.fu-berlin.de}.
$^1$ Mathematical Institute of the Hungarian Academy of Sciences, 
Budapest, Hungary. e-mail: {\tt soukup@math-inst.hu}. The second author 
was supported by the Hungarian National Foundation for Scientific 
Research grant no.\ 16391.}

\section{Introduction}

For a regular $\kappa$, 
a \po\ $P$ is said to have the $\kappa$ Freese-Nation property (the 
$\kappa$-FN for short) if there is a mapping ($\kappa$-FN mapping) 
$\mapping{f}{P}{[P]^{<\kappa}}$ \st\ for any $p$, $q\in P$ if 
$p\leq q$ then there is $r\in f(p)\cap f(q)$ \st\ $p\leq r\leq q$. 

Freese and Nation \cite{freese-nation} used the $\aleph_0$-FN 
in a characterization of projective lattices and asked if this property 
alone already characterizes the projectiveness. L.\ 
Heindorf gave a 
negative answer to the question showing that the \Bas\ 
with the $\aleph_0$-FN are exactly those which are openly generated. It 
is known that the class of openly generated \Bas\ contains projective 
\Bas\ as a proper subclass (see \cite{heindorf-shapiro} --- openly 
generated \Bas\ are called `rc-filtered' there). Heindorf and Shapiro 
\cite{heindorf-shapiro} then studied the $\aleph_1$-FN which they called 
the weak Freese-Nation property and proved some elementary properties of 
the \Bas\ with this property. \Pos\ with the $\kappa$-FN for arbitrary 
regular $\kappa$ were further studied in Fuchino, Koppelberg and Shelah 
\cite{WFN}. Koppelberg \cite{koppelberg} gives some nice applications of 
the $\aleph_1$-FN. 

In the following we shall quote some elementary facts from \cite{WFN} 
which we need later. First of all, it can be readily seen that every 
small \po\ has the $\kappa$-FN:
\begin{Lemma}{\rm (\cite{WFN})} \label{baka}
Every \po\ $P$ of cardinality $\leq\kappa$ has the $\kappa$-FN.\qed
\end{Lemma}

For a \po\ $P$ and a sub-ordering $Q\subseteq P$, we say that $Q$ is a 
$\kappa$-sub-ordering of $P$ and denote it with $Q\leq_\kappa P$ if, for 
every $p\in P$, the set $\setof{q\in Q}{q\leq p}$ has a cofinal subset of 
cardinality $<\kappa$ and the set $\setof{q\in Q}{q\geq p}$ has a 
coinitial subset of cardinality $<\kappa$. 
\begin{Lemma}{\rm (\cite{WFN})} \label{conti-chain}Suppose that 
$\delta$ is a limit ordinal and $(P_\alpha)_{\alpha\leq\delta}$ a 
continuously increasing chain of \pos\ \st\ 
$P_\alpha\leq_\kappa P_\delta$ for all $\alpha<\delta$. 
If $P_\alpha$ has the $\kappa$-FN  
for every $\alpha<\delta$, then 
$P_\delta$ also has the $\kappa$-FN.\qed
\end{Lemma}
For application of \Lemmaof{conti-chain}, it is enough to have 
$P_\alpha\leq_\kappa P_\delta$ and the $\kappa$-FN of $P_\alpha$ for every 
$\alpha<\delta$ \st\ either $\alpha$ is a successor or of cofinality 
$\geq\kappa$: $P_\alpha\leq_\kappa P_\delta$ for $\alpha<\delta$ of 
cofinality 
$<\kappa$ follows from this since such $P_\alpha$ can be represented as 
the union of $<\kappa$ many $\kappa$-sub-orderings of $P_\delta$. Hence by 
inductive application of \Lemmaof{conti-chain}, we can show that 
$P_\alpha$ satisfies the $\kappa$-FN for every $\alpha\leq\delta$. 
Similarly, if $\delta$ is a cardinal $>\kappa$, then it is enough to have 
$P_\alpha<_\kappa P_\delta$ and the $\kappa$-FN of $P_\alpha$ for every limit 
$\alpha<\delta$ of cofinality $\geq\kappa$. 
\begin{Prop}{\rm(\cite{WFN})} \label{characterization}
For a regular $\kappa$ and a \po\ $P$, 
\tfae:\smallskip\\
\assert{1} $P$ has the $\kappa$-FN;\smallskip\\
\assert{2} For some, or equivalently, any sufficiently large $\chi$, if 
$M\prec\calH_\chi=(\calH_\chi,\in)$ is \st\ 
$P\in M$, $\kappa\subseteq M$ and $\cardof{M}=\kappa$ then 
$P\cap M\leq_\kappa P$ 
holds;\smallskip\\
\assert{3} $\setof{C\in[P]^{\kappa}}{C\leq_\kappa P}$ contains a club 
set.\qed 
\end{Prop}
Though \Propof{characterization},\assertof{2} is quite useful to show that 
a \po\ has the $\kappa$-FN, sometimes it is quite difficult to 
check \Propof{characterization},\,\assertof{2} as in the case of 
the $\aleph_1$-FN of $P(\omega)$ or $[\kappa]^{<\omega}$: in these cases 
it is independent if \Propof{characterization},\,\assertof{2} holds. 
Applications like \Corof{appl-of-char} in mind, 
we could think of 
another possible variant of \Propof{characterization},\,\assertof{2} in 
terms of the following 
weakening of the notion of internally approachability from \cite{MM}:
for a regular $\kappa$ and a sufficiently large $\chi$, we shall call an 
elementary submodel $M$ of 
$\calH_\chi$ {\em $V_{\kappa}$-like} if, either 
$\kappa=\aleph_0$ and $M$ is countable, or there is an 
increasing 
sequence $(M_\alpha)_{\alpha<\kappa}$ of elementary submodels 
of $M$ of cardinality less than $\kappa$ \st\ 
$M_\alpha\in M_{\alpha+1}$ 
for all 
$\alpha<\kappa$ and $M=\bigcup_{\alpha<\kappa}M_\alpha$. 

In \cite{WFN}, a characterization of the $\kappa$-FN using 
$V_\kappa$-like elementary submodels in place of elementary submodels in 
\Propof{characterization},\,\assertof{2} was discussed. Unfortunately it 
appeared that some consequences of $\neg0^\#$ are necessary for the 
characterization (see ``Added in Proof'' in \cite{WFN}). 
In this paper, we introduce a weakening of the very weak square principle 
from \cite{foreman-magidor} --- the principle $\Box^{***}_{\kappa,\mu}$. 
In section \ref{jensen-matrix} we show the equivalence of 
$\Box^{***}_{\kappa,\mu}$ with the existence of a matrix 
$(M_{\alpha,\beta})_{\alpha<\mu^+,\beta<\cf(\mu)}$ 
--- which we called (weak) $(\kappa,\mu)$-Jensen matrix --- of elementary 
submodels of $\calH(\chi)$ for sufficiently large $\chi$ with certain 
properties. This fact is used in section \ref{xchar} to show 
that $\Box^{***}_{\kappa,\lambda}$ together with a very weak version of 
the Singular Cardinals Hypothesis yields the characterization of \pos\ 
with the $\kappa$-FN in terms of 
$V_\kappa$-like elementary 
submodels (\Thmof{correction}). ZFC or even ZFC $+$ GCH is not enough for 
this characterization: in section \ref{chang}, we show that, under Chang's 
conjecture for $\aleph_\omega$, there is a counter-example to the 
characterization. Together with \Thmof{correction}, this counter-example also 
shows that $\Box^{***}_{\aleph_1,\aleph_\omega}$ is not a theorem in ZFC 
$+$ GCH. 

One of the most natural questions concerning the $\kappa$-FN would be if 
$(\powersetof{\omega},\subseteq)$ has the $\aleph_1$-FN. It is easy to see 
that $(\powersetof{\omega},\subseteq)$ has the $\aleph_1$-FN iff 
$(\powersetof{\omega}/fin,\subseteq^*)$ does (see \cite{WFN}). 
See also Koppelberg \cite{koppelberg} for some 
consequences of the $\aleph_1$-FN of $\powersetof{\omega}/fin$. By 
\Lemmaof{baka}, $\powersetof{\omega}$ has the $\aleph_1$-FN under CH. 
In section \ref{cohen}, we show that $\powersetof{\omega}$ and a lot of 
other ccc \cBas\ still have the $\aleph_1$-FN in a model obtained by 
adding arbitrary number of Cohen reals to a model of, say, $V=L$. On the 
other hand it can happen very easily that $\powersetof{\omega}$ does not 
have the $\aleph_1$-FN. In \cite{WFN}, it was shown that this is the case 
when ${\bf b}>\aleph_1$ or, more generally, if there is a 
$\subseteq^*$-sequence of elements of $\psof{\omega}$ of order-type 
$\geq\omega_2$. In section \ref{lusin}, we show that the 
existence of an $\aleph_2$-Lusin gap can be another reason for failure 
of the $\aleph_1$-FN. At the moment the authors do not know if there are 
yet other reasons for failure of the $\aleph_1$-FN of $\powersetof{\omega}$:
\begin{Problem}
Suppose that $\powersetof{\omega}$ does not have any increasing chain of 
length $\geq\aleph_2$ \wrt\ $\subseteq^*$ and that there is no 
$\aleph_2$-Lusin gap. Does it follow that $\powersetof{\omega}$ has the 
$\aleph_1$-FN\,? 
\end{Problem}
\par
Our notation is fairly standard. The following are possible deviations 
from the standard: for $C\subseteq\kappa$, we denote with 
$(C)'$ the set of limit points of $C$ other than $\kappa$. For an ordinal 
$\alpha$, 
$\Lim(\alpha)=\setof{\beta<\alpha}{\beta\xmbox{ is a limit ordinal}}$. 
For a \po\ 
$P$, 
$\cf(P)=\min\setof{\cardof{X}}{X\subseteq P,X\mbox{ is cofinal in }P}$.
$[\lambda]^{<\kappa}=\setof{X\subseteq\lambda}{\cardof{X}<\kappa}$ is 
often seen as the \po\ $([\lambda]^{<\kappa},\subseteq)$. 
If $Q$ is a sub-ordering of a \po\ $P$ and $p\in P$ then 
$Q\upper p=\setof{q\in Q}{q\geq p}$ and $Q\restr p=\setof{q\in Q}{q\leq p}$. 

\section{Very weak square and Jensen matrix}\label{jensen-matrix}
For a cardinal $\mu$, the {\em weak square principle for $\mu$ }
(notation: $\Box^*_\mu$) 
is the statement: there is a sequence 
$(\calC_\alpha)_{\alpha\in\Lim(\mu^+)}$ \st\ for every 
$\alpha\in\Lim(\mu^+)$
\begin{assertion}{}\mbox{}%
\lassert{w1} $\calC_\alpha\subseteq\powersetof{\alpha}$ and 
$\cardof{\calC_\alpha}\leq\mu$;\smallskip\\
\lassert{w2} every $C\in\calC_\alpha$ is club in $\alpha$ and if 
$\cf(\alpha)<\mu$ then $\otp(C)<\mu$;\smallskip\\
\lassert{w3} there is $C\in\calC_\alpha$ \st\ 
for every $\delta\in(C)'$, $C\cap\delta\in\calC_\delta$.
\end{assertion}
Clearly we have $\Box_\mu\rightarrow\,\Box^*_\mu$. Jensen \cite{jensen} 
proved that $\Box^*_\mu$ is equivalent to the existence of a special 
Aronszajn tree on $\mu^+$. Ben-David and Magidor \cite{ben-david-magidor} 
showed that the weak square principle for a singular $\mu$ is actually 
weaker than the square 
principle: they constructed a model of $\Box^*_{\aleph_\omega}$ and 
$\neg\Box_{\aleph_\omega}$ starting from a model with a supercompact 
cardinal. 

Foreman and Magidor considered in \cite{foreman-magidor} the following 
principle which is, e.g.\ under GCH, a weakening of $\Box^*$ principle: 
for a cardinal $\mu$, {\em the very weak square principle for 
$\mu$} holds if there is a sequence\ 
$(C_\alpha)_{\alpha<\mu^+}$ and a 
club $D\subseteq\mu^+$ \st\ for every $\alpha\in D$
\begin{assertion}{}\mbox{}%
\lassert{v1} $C_\alpha\subseteq\alpha$, $C_\alpha$ is unbounded in 
$\alpha$;\smallskip\\ 
\lassert{v2} for all bounded $x\in[C_\alpha]^{<\omega_1}$, there is 
$\beta<\alpha$ \st\ $x=C_\beta$.
\end{assertion}
In this paper, 
we shall use the following yet weaker variant of the very weak 
square principle. 
For a regular cardinal $\kappa$ and $\mu>\kappa$, let 
$\Box^{***}_{\kappa,\mu}$ be the following assertion: there exists a 
sequence $(C_\alpha)_{\alpha<\mu^+}$ and a club set 
$D\subseteq\mu^+$ \st\ for $\alpha\in D$ with $\cf(\alpha)\geq\kappa$
\begin{assertion}{}\mbox{}%
\lassert{y1} $C_\alpha\subseteq\alpha$, $C_\alpha$ is unbounded in 
$\alpha$;\smallskip\\ 
\lassert{y2} $[\alpha]^{<\kappa}\cap\setof{C_{\alpha'}}{\alpha'<\alpha}$ 
dominates $[C_\alpha]^{<\kappa}$ (\wrt\ $\subseteq$). 
\end{assertion}
Since \assertof{y2} remains valid when $C_\alpha$'s for $\alpha\in D$ are 
shrunk, we may replace \assertof{y1} by
\begin{assertion}{}\mbox{}%
\lassert{y1'} $C_\alpha\subseteq\alpha$, $C_\alpha$ is unbounded in 
$\alpha$ and $\otp(C_\alpha)=\cf(\alpha)$.
\end{assertion}
A corresponding 
remark holds also for the sequence of the very weak square principle. 
\begin{Lemma}
\assert{a} The very weak square principle for $\mu$ implies 
$\Box^{***}_{\omega_1,\mu}$. \smallskip\\
\assert{b} For a cardinal $\mu$ and a regular $\kappa$ \st\ 
$\cf(\mu)<\kappa$ and 
$\cf([\lambda]^{<\kappa},\subseteq)\leq\mu$ for every $\lambda<\mu$, 
$\Box^*_\mu$ implies $\Box^{***}_{\kappa,\mu}$.
\end{Lemma}
\prf \assertof{a} is clear. For \assertof{b}, let 
$(\calC_\alpha)_{\alpha\in\Lim(\mu^+)}$ be a weak square sequence. Let 
$\calC_\alpha=\setof{C_{\alpha,\beta}}{\beta<\mu}$ for every 
$\alpha\in\Lim(\mu^+)$. \Wolog, we may assume that $C_{\alpha,0}$ is as 
the $C$ 
in \assertof{w3}. By shrinking $C_{\alpha,\beta}$'s if necessary, we may 
also assume that 
$\cardof{C_{\alpha,\beta}}<\mu$ for every $\alpha\in\Lim(\mu^+)$ and 
$\beta<\mu$. 
For $\alpha\in\Lim(\mu^+)$ and $\beta<\mu$, 
let 
$X_{\alpha,\beta}$ be a cofinal subset of $[C_{\alpha,\beta}]^{<\kappa}$ 
of cardinality $\leq\mu$ and let 
$\setof{C_\alpha}{\alpha\in\mu^+\setminus\Lim(\mu^+)}$ be an enumeration 
of $\bigcup\setof{X_{\alpha,\beta}}{\alpha\in\Lim(\mu^+),\beta<\mu}$. For each 
$\alpha\in\Lim(\mu^+)$, let $C_\alpha=C_{\alpha,0}$. Let 
$\mapping{F}{\mu^+}{\mu^+}$ be defined by 
$F(\xi)=
	\min\setof{\gamma<\mu^+}{%
		X_{\xi,\beta}\subseteq\setof{C_\alpha}{\alpha<\gamma}
		\xmbox{ for every }\beta<\mu}$.
Let $D\subseteq\mu^+$ be a club set closed \wrt\ $F$. Then 
$(C_\alpha)_{\alpha<\mu^+}$ and $D$ are as in the definition of 
$\Box^{***}_{\kappa,\mu}$. To see that \assertof{y2} is satisfied, let 
$\alpha\in D$ be \st\ $\cf(\alpha)\geq\kappa$ and 
$x\in[C_\alpha]^{<\kappa}$. By definition we have 
$C_\alpha=C_{\alpha,0}$. Hence there are 
$\alpha'\in\alpha\cap\Lim(\mu^+)$ and 
$\beta<\lambda$ \st\ $x\in[C_{\alpha',\beta}]^{<\kappa}$. Since 
$\alpha$ is closed \wrt\ $F$, there is some $\gamma<\alpha$ \st\ 
$x\subseteq C_\gamma\in [C_{\alpha',\beta}]^{<\kappa}$. 	
\qedofLemma
\qedskip\par
$\Box^{***}_{\kappa,\mu}$ has some 
influence on cardinal arithmetic of cardinals below $\mu$:
\begin{Lemma}\label{card-arith}
Suppose that $\kappa$ is regular and $\mu$ is \st\ $\cf(\mu)<\kappa$. 
If $\Box^{***}_{\kappa,\mu}$ holds, 
then we have $\cf([\lambda]^{<\kappa},\subseteq)<\mu$ for every 
$\lambda<\mu$. 
\end{Lemma}
\prf
Let $(C_\alpha)_{\alpha<\mu^+}$, and $D\subseteq\mu^+$ be witnesses of 
$\Box^{***}_{\kappa,\mu}$. For $\lambda<\mu$, let $\delta\in D$ be \st\ 
$\cf(\delta)\geq\lambda+\kappa$. Then 
$\setof{C_\alpha}{\alpha<\delta}\cap[\delta]^{<\kappa}$ is cofinal in 
$[C_\delta]^{<\kappa}$. Since $\cardof{\delta}\leq\mu$ it follows that 
$\cf([C_\delta]^{<\kappa},\subseteq)\leq\mu$. As the order type of 
$C_\delta$ is at least $\lambda$, it follows that 
$\nu=\cf([\lambda]^{<\kappa},\subseteq)\leq\mu$. But 
$\cf(\nu)\geq\kappa$. Hence $\nu<\mu$. 
\qedofLemma\qedskip\\

Suppose now that $\kappa$ is a regular cardinal and $\mu>\kappa$ is \st\ 
$\cf(\mu)<\kappa$. Let $\mu^*=\cf(\mu)$. For a sufficiently large $\chi$ 
and $x\in\calH(\chi)$, let us call a sequence
$(M_{\alpha,\beta})_{\alpha<\mu^+,\beta<\mu^*}$ a {\em$(\kappa,\mu)$-Jensen 
matrix over $x$} --- or just {\em Jensen matrix over $x$} if it is clear from 
the context which $\kappa$ and $\mu$ are meant --- if the following 
conditions hold:
\begin{assertion}{}\mbox{}%
\lassert{j1} $M_{\alpha,\beta}\prec\calH(\chi)$, 
$x\in M_{\alpha,\beta}$, $\kappa+1\subseteq M_{\alpha,\beta}$ and 
$\cardof{M_{\alpha,\beta}}<\mu$ for all $\alpha<\mu^+$ and 
$\beta<\mu^*$\,;\smallskip\\
\lassert{j2} $(M_{\alpha,\beta})_{\beta<\mu^*}$ is an increasing sequence 
for each $\alpha<\mu^+$\,;\smallskip\\
\lassert{j3} if $\alpha<\mu^+$ is \st\ $\cf(\alpha)\geq\kappa$, then there 
is $\beta^*<\mu^*$ \st, for every $\beta^*\leq\beta<\mu^*$, 
$[M_{\alpha,\beta}]^{<\kappa}\cap M_{\alpha,\beta}$ is cofinal in 
$([M_{\alpha,\beta}]^{<\kappa},\subseteq)$\,;
\end{assertion}
For $\alpha<\mu^+$, let
$M_\alpha=\bigcup_{\beta<\mu^*}M_{\alpha,\beta}$. By \assertof{j1} and 
\assertof{j2}, we have 
$M_\alpha\prec\calH(\chi)$. 
\begin{assertion}{}\mbox{}%
\lassert{j4} $(M_\alpha)_{\alpha<\mu^+}$ is continuously increasing and $\mu^+\subseteq\bigcup_{\alpha<\mu^+}M_\alpha$.
\end{assertion}
The choice of the term ``Jensen matrix'' was suggested by 
\cite{foreman-magidor} in which (under GCH) a matrix of subsets of $\mu^+$ 
having some properties similar to those of the sequence 
$(\mu^+\cap M_{\alpha,\beta})_{\alpha<\mu^+,\beta<\mu^*}$ for 
$(M_{\alpha,\beta})_{\alpha<\mu^+,\beta<\mu^*}$ as above is called a 
Jensen sequence. 
Note that, in the case of $2^{<\kappa}=\kappa$, \assertof{j3} can be 
replaced by the following seemingly stronger property:
\begin{assertion}{}\mbox{}%
\lassert{j3'} if $\alpha<\mu^+$ is \st\ $\cf(\alpha)\geq\kappa$, then there 
is $\beta^*<\mu^*$ \st, for every $\beta^*\leq\beta<\mu^*$, 
$[M_{\alpha,\beta}]^{<\kappa}\subseteq M_{\alpha,\beta}$.
\end{assertion}
This is simply because of the following observation:
\begin{Lemma}\label{easy}
Suppose that $2^{<\kappa}=\kappa$ and $M$ is an elementary submodel of 
$\calH(\chi)$ for some sufficiently large $\chi$ and 
$\kappa\subseteq M$. If $[M]^{<\kappa}\cap M$ is cofinal in 
$[M]^{<\kappa}$, then we have $[M]^{<\kappa}\subseteq M$.
\end{Lemma}
\prf
Let $x\in[M]^{<\kappa}$. We show that $x\in M$. By assumption there is 
$y\in [M]^{<\kappa}\cap M$ \st\ $x\subseteq y$. Let 
$\eta=\cardof{\psof{y}}$. Then $\eta\in M$ and $\eta\leq\kappa$. Let 
$(y_\alpha)_{\alpha<\eta}\in M$ be an enumeration of $\psof{y}$. Then 
there is an $\alpha_0<\eta$ \st\ $x=y_{\alpha_0}$. But since 
$\kappa\subseteq M$, we have $\alpha_0\in M$ and $y_{\alpha_0}\in M$ as well. 
\qedofLemma\qedskip
\\
Note also that, if $M\prec\calH(\chi)$ is $V_\kappa$-like, then 
$[M]^{<\kappa}\cap M$ is cofinal in $[M]^{<\kappa}$. Hence, under 
$2^{<\kappa}=\kappa$, $M\prec\calH(\chi)$ is $V_\kappa$-like if and only 
if $\cardof{M}=\kappa$ and $[M]^{<\kappa}\subseteq M$.

In the following theorem, we show that $\Box^{***}_{\kappa,\mu}$ together 
with a very 
weak version of the Singular Cardinals Hypothesis below $\mu$ implies the 
existence of a Jensen matrix: 
\begin{Thm}\label{Jensen-exists}
Suppose that $\kappa$ is a regular cardinal and $\mu>\kappa$ is \st\ 
$\cf(\mu)<\kappa$. If we have $\cf([\lambda]^{<\kappa},\subseteq)=\lambda$ 
for cofinally many 
$\lambda<\mu$ and $\Box^{***}_{\kappa,\mu}$ holds, then, 
for any sufficiently large $\chi$ and $x\in\calH(\chi)$, there is a 
$(\kappa,\mu)$-Jensen matrix over $x$.
\end{Thm}
\prf
Let $\mu^*=\cf(\mu)$ and $(\mu_\beta)_{\beta<\mu^*}$ be an increasing 
sequence of cardinals below $\mu$ \st\ $\mu_0>\mu^*$, 
$\sup\setof{\mu_\beta}{\beta<\mu^*}=\mu$ 
and $\cf([\mu_\beta]^{<\kappa},\subseteq)=\mu_\beta$ for every 
$\beta<\mu^*$. 
Let $(C_\alpha)_{\alpha\in\mu^+}$ and $D\subseteq\mu^+$ be as in the 
definition of $\Box^{***}_{\kappa,\mu}$. 
\Wolog, we may assume that $\cardof{C_\alpha}\leq\cf(\alpha)$ for all 
$\alpha>\mu^+$. We may also assume that $\alpha>\mu$ for every $\alpha\in D$.

In the following, we fix a well ordering $\unlhd$ on $\calH(\chi)$ and, 
when we talk about $\calH(\chi)$ as a structure, we mean 
$\calH(\chi)=(\calH(\chi),\in,\unlhd)$. 
$X\subseteq\calH(\chi)$ as a substructure of $\calH(\chi)$ is thus the 
structure $(X,\in\cap X^2,\unlhd\cap X^2)$ --- for notational simplicity 
we shall denote such a structure simply by $(X,\in,\unlhd)$.

Let $N\in\calH(\chi)$ be an elementary substructure of 
$\calH(\chi)$ \st\ $N$ contains every 
thing needed below --- in particular, we let $\mu^+\subseteq N$ and 
$x$, $(C_{\alpha})_{\alpha<\mu^+}$, $D$, 
$(\mu_\alpha)_{\alpha<\mu^*}\in N$. Let $(N_\xi)_{\xi<\kappa}$ be an 
increasing sequence of elementary submodels of $\calH(\chi)$ \st\ 
\begin{assertion}{}\mbox{}%
\lassert{0} $N_0=N$\,; \\
\lassert{1} $N_\xi\in \calH(\chi)$ for every $\xi<\kappa$ and \\
\lassert{2} $(N_\eta)_{\eta\leq\xi}\in N_{\xi+1}$ for every $\xi<\kappa$.
\end{assertion}

Now, for each $\xi<\kappa$, let 
\[ \calN_\xi
	=(N_\xi,\in,\unlhd,R_\xi,\mu,\mu^*,
		D,
		(C_{\alpha})_{\alpha<\mu^+},
		(\mu_\alpha)_{\alpha<\mu^*},x,\eta)_{\eta<\xi}
\]\noindent
where $R_\xi$ is the relation $\setof{(\eta, N_\eta)}{\eta<\xi}$. 
For $X\subseteq\mu^+$, let us denote with $sk_\xi(X)$ the Skolem-hull of 
$X$ \wrt\ the built-in Skolem functions of $\calN_\xi$ (induced from 
$\unlhd$). For $\xi<\xi'<\kappa$, $\calN_\xi$ is an element of 
$\calN_{\xi'}$ by \assertof{2} and the Skolem functions of $\calN_\xi$ 
are also 
elements of $\calN_{\xi'}$. In particular, we have 
$sk_\xi(X)\subseteq sk_{\xi'}(X)$. It follows that  
$sk(X)=\bigcup_{\xi<\kappa}sk_\xi(X)$ is an elementary submodel of 
$\calH(\chi)$. Note also that, if $X\subseteq\mu^+$ is an element of 
$sk_{\xi'}(Y)$ then, since $sk_\xi(X)$ is definable in $sk_{\xi'}(Y)$, we 
have $sk_\xi(X)\in sk_{\xi'}(Y)$. 

For the proof of the theorem, 
it is clearly enough to construct $M_{\alpha,\beta}$ with \assertof{j1} 
--- \assertof{j4} for every $\alpha$ in the club set 
$D$ and for every 
$\beta<\mu^*$. Let
\[ M_{\alpha,\beta}=sk(\mu_\beta\cup C_{\alpha})
\]\noindent
for $\alpha\in D$ and $\beta<\mu^*$. We show that 
$(M_{\alpha,\beta})_{\alpha\in D,\beta<\mu^*}$ is as desired. It is clear 
that \assertof{j1} and \assertof{j2} hold. We need the following claim to 
show the other properties: 
\begin{Claim}\label{claim-alpha}
$M_\alpha=sk(\alpha)$ for every $\alpha\in D$.
\end{Claim}
\prfofClaim
``$\subseteq$'' is clear since $\mu_\beta\cup C_{\alpha}\subseteq\alpha$ 
for every $\alpha\in D$ and $\beta<\mu^*$. For ``$\supseteq$'', it is 
enough to show that $\alpha\subseteq M_\alpha$. Let 
$\gamma<\alpha$. By \assertof{y1}, there is $\gamma_1\in C_{\alpha}$ 
\st\ $\gamma<\gamma_1$. Let $f\in M_{\alpha,0}$ be a surjection from 
$\mu$ to $\gamma_1$, and let $\delta<\mu$ be \st\ $f(\delta)=\gamma$. 
Then we have 
$\gamma=f(\delta)\in M_{\alpha,\beta^*}\subseteq M_\alpha$ for 
$\beta^*<\mu^*$ \st\ $\delta<\mu_{\beta^*}$. 
\qedofClaim
\smallskip

For \assertof{j3}, suppose that $\alpha\in D$ and 
$\cf(\alpha)\geq\kappa$. Let $\beta^*$ be \st\ 
$\cardof{C_\alpha}<\mu_{\beta^*}$ and 
$[\alpha]^{<\kappa}\cap
	\setof{C_{\alpha'}}{\alpha'<\alpha,\,\alpha'\in M_{\alpha,\beta^*}}$ 
dominates 
$[C_\alpha]^{<\kappa}$. The last property is possible by \assertof{y2} 
and \Claimabove. 
We show that 
this $\beta^*$ is as needed in \assertof{j3}. Let $\beta<\mu^*$ be \st\ 
$\beta^*\leq\beta$ and 
suppose that $x\in[M_{\alpha,\beta}]^{<\kappa}$. Then there are 
$u\in[\mu_\beta]^{<\kappa}$ and $v\in[C_{\alpha}]^{<\kappa}$ \st\ 
$x\subseteq sk(u\cup v)$. Since $\mu_\beta\in M_{\alpha,\beta}$ and 
$\cf([\mu_\beta]^{<\kappa},\subseteq)=\mu_\beta$, there is 
$X\in M_{\alpha,\beta}$ \st\ $X\subseteq[\mu_\beta]^{<\kappa}$, 
$\cardof{X}=\mu_\beta$ and $X$ is cofinal in 
$([\mu_\beta]^{<\kappa},\subseteq)$. Since 
$\mu_\beta\subseteq M_{\alpha,\beta}$, it follows that 
$X\subseteq M_{\alpha,\beta}$. Hence there is 
$u'\in M_{\alpha,\beta}\cap[\mu_\beta]^{<\kappa}$ 
\st\ $u\subseteq u'$. On the other hand, by definition of $\beta^*$, 
$\cf(\alpha)\geq\kappa$, there is 
$\alpha'\in \alpha\cap M_{\alpha,\beta}$ \st\ 
$C_{\alpha'}\in [\alpha]^{<\kappa}$ and $v\subseteq C_{\alpha'}$. We have 
$x\subseteq sk(u\cup v)\subseteq sk(u'\cup C_{\alpha'})$. By regularity of 
$\kappa$, there is $\xi<\kappa$ \st\  
$x\subseteq sk_\xi(u'\cup C_{\alpha'})$. But 
$sk_\xi(u'\cup C_{\alpha'})\in M_{\alpha,\beta}$ and 
$\cardof{sk_\xi(u'\cup C_{\alpha'})}<\kappa$.

\assertof{j4} follows immediately from \Claimof{claim-alpha}.
\qedofThm 
\qedskip\\
Note that, in the proof above, the sequence 
$(M_{\alpha,\beta})_{\alpha<\mu^+,\beta<\mu^*}$ satisfies also: 
\begin{assertion}{}\mbox{}%
\lassert{j5} for $\alpha<\alpha'<\mu^+$ and $\beta<\mu^*$, there is 
$\beta'<\mu^*$ \st\ 
$M_{\alpha,\beta}\subseteq M_{\alpha',\beta'}$.
\end{assertion}
[\,Suppose that $\alpha$, $\alpha'\in D$ are \st\ 
$\alpha<\alpha'$ and $\beta<\mu^*$. By \Claimof{claim-alpha}, there is 
$\beta'<\mu^*$ \st\ $\beta<\beta'$, $\alpha\in M_{\alpha',\beta'}$
and $\otp(C_{\alpha})\leq\mu_{\beta'}$. 
Then we have 
$C_{\alpha}\in M_{\alpha',\beta'}$ and 
$C_{\alpha}\subseteq M_{\alpha',\beta'}$. Also 
$\mu_\beta\in M_{\alpha',\beta'}$ and 
$\mu_\beta\subseteq\mu_{\beta'}\subseteq M_{\alpha',\beta'}$. Hence 
it follows that 
$M_{\alpha,\beta}=sk(\mu_\beta\cup C_{\alpha})
	\subseteq M_{\alpha',\beta'}$\,.]

Conversely, the existence of a $(\kappa,\mu)$-Jensen matrix (over 
some/any $x$) implies $\Box^{***}_{\kappa,\mu}$:
\begin{Thm}\label{converse}
Suppose that $\kappa$ is a regular cardinal and $\mu>\kappa$ is \st\ 
$\cf(\mu)<\kappa$. If there exists a $(\kappa,\mu)$-Jensen matrix, then 
$\Box^{***}_{\kappa,\mu}$ holds.
\end{Thm}
\prf
Let $\mu^*=\cf(\mu)$ and $(M_{\alpha,\beta})_{\alpha<\mu^+,\beta<\mu^*}$ 
be a $(\kappa,\mu)$-Jensen matrix. Let 
$M_\alpha=\bigcup_{\beta<\mu^*}M_{\alpha,\beta}$ for each $\alpha<\mu^+$. For 
\[ \textstyle X=
	\bigcup\setof{[M_{\alpha,\beta}]^{<\kappa}\cap M_{\alpha,\beta}}{%
		\alpha<\mu^+,\,\beta<\mu^*},
\]\noindent
let $\setof{C_{\alpha+1}}{\alpha<\mu^+}$ be an enumeration of $X$. Let 
\[ D=
	\setof{\alpha<\mu^+}{%
		\begin{array}[t]{@{}l}
			M_\alpha\cap\mu^+=\alpha,\\
			\setof{C_{\alpha'+1}}{\alpha'<\alpha}
				\supseteq 
				[M_{\alpha^\dagger,\beta}]^{<\kappa}\cap 
					M_{\alpha^\dagger,\beta}\\
			\mbox{for every }\alpha^\dagger<\alpha,\,\beta<\mu^*\,\,}.
		\end{array}
\]\noindent
By \assertof{j1} and \assertof{j4}, $D$ is a club subset of $\mu^+$. For 
$\alpha\in D$ with $\cf(\alpha)\geq\kappa$, let 
$C_\alpha=M_{\alpha,\beta_\alpha}\cap \alpha$ where $\beta_\alpha<\mu^+$ be 
\st\ $M_{\alpha,\beta_\alpha}\cap\alpha$ is cofinal in $\alpha$ (this is 
possible as $M_\alpha\cap \mu^+=\alpha$ and $\mu^*<\kappa$) and that 
$[M_{\alpha,\beta_\alpha}]^{<\kappa}\cap M_{\alpha,\beta_\alpha}$ is 
cofinal in $[M_{\alpha,\beta_\alpha}]^{<\kappa}$ (possible by 
\assertof{j3}\,). For 
$\alpha\in\Lim(\mu^+)\setminus\setof{\alpha\in D}{\cf(\alpha)\geq\kappa}$, 
let $C_\alpha=\emptyset$. We claim that $(C_\alpha)_{\alpha<\mu^+}$ and 
$D$ as above satisfy the conditions in the definition of 
$\Box^{***}_{\kappa,\mu}$: \assertof{y1} is clear by definition of 
$C_\alpha$'s. To show \assertof{y2}, let $\alpha\in D$ be \st\ 
$\cf(\alpha)\geq\kappa$ and $x\in[C_\alpha]^{<\kappa}$. Then, by 
the choice of $\beta_\alpha$, there is 
$y\in[\alpha]^{<\kappa}\cap M_{\alpha,\beta_\alpha}$ \st\ 
$x\subseteq y$. By \assertof{j4}, there are $\alpha^\dagger<\alpha$ and 
$\beta^\dagger<\mu^*$ \st\ $y\in M_{\alpha^\dagger,\beta^\dagger}$. By 
definition of 
$D$, it follows that $y=C_{\alpha'+1}$ for some $\alpha'<\alpha$. 
Thus $[\alpha]^{<\kappa}\cap\setof{C_{\alpha'}}{\alpha'<\alpha}$ 
dominates $[C_\alpha]^{<\kappa}$. 
\qedofThm\qedskip
\\
Thus, if $\cf([\lambda]^{<\kappa},\subseteq)=\lambda$ for cofinally many 
$\lambda<\mu$, $\Box^{***}_{\kappa,\mu}$ is equivalent to the existence 
of a $(\kappa,\mu)$-Jensen matrix. Using the following weakening of the 
notion of Jensen matrix, we can obtain a characterization of 
$\Box^{***}_{\kappa,\mu}$ in ZFC: for a regular cardinal $\kappa$ and 
$\mu>\kappa$ \st\ $\mu^*=\cf(\mu)<\kappa$, let us call a matrix 
$(M_{\alpha,\beta})_{\alpha<\mu^+,\beta<\mu^*}$ of elementary submodels 
of $\calH(\chi)$ for a sufficiently large $\chi$, {\em a weak 
$(\kappa,\mu)$-Jensen matrix over $x$}, if it satisfies \assertof{j1}, 
\assertof{j2}, \assertof{j4} for 
$M_\alpha=\bigcup_{\beta<\mu^*}M_{\alpha,\beta}$, $\alpha<\mu^+$, and
\begin{assertion}{}\mbox{}%
\lassert{j3$^-$\relax} if $\alpha<\mu^+$ is \st\ 
$\cf(\alpha)\geq\kappa$, then there 
is $\beta^*<\mu^*$ \st, for every $\beta^*\leq\beta<\mu^*$, 
$[M_{\alpha,\beta}]^{<\kappa}\cap M_{\alpha}$ is cofinal in 
$([M_{\alpha,\beta}]^{<\kappa},\subseteq)$\,.
\end{assertion}
Since $\mu^*<\kappa$, the last condition is equivalent with:
\begin{assertion}{}\mbox{}%
\lassert{j3$^*$\relax} if $\alpha<\mu^+$ is \st\ 
$\cf(\alpha)\geq\kappa$, then there 
is $\beta^*<\mu^*$ \st, for every $\beta^*\leq\beta<\mu^*$, 
there is $\beta'<\mu^*$ \st\ 
$[M_{\alpha,\beta}]^{<\kappa}\cap M_{\alpha,\beta'}$ is cofinal in 
$([M_{\alpha,\beta}]^{<\kappa},\subseteq)$\,.
\end{assertion}
\begin{Thm}
Suppose that $\kappa$ is a regular cardinal and $\mu>\kappa$ is \st\ 
$\mu^*=\cf(\mu)<\kappa$. Then $\Box^{***}_{\kappa,\mu}$ holds if and only if 
there is a weak $(\kappa,\mu)$-Jensen matrix over some/any $x$. 
\end{Thm}
\prf
For the forward direction the proof is almost the same as that of 
\Thmof{Jensen-exists}. We let $(\mu_\beta)_{\beta<\mu^*}$ here merely an 
increasing sequence of regular cardinals with the limit $\mu$. Then 
$(M_{\alpha,\beta})_{\alpha\in D,\beta<\mu^*}$ is constructed just as in 
the proof of \Thmof{Jensen-exists}. \Lemmaof{card-arith} is then used 
to see that \assertof{j3$^-$\relax} is satisfied by this matrix. For the 
converse, just the same proof as that of \Thmof{converse} will do.
\qedofThm\qedskip\\
Existence of a Jensen-matrix is not a theorem in ZFC: we show in section 
\ref{chang} that 
the Chang's conjecture for $\aleph_\omega$ together with 
$2^{\aleph_\omega}=\aleph_{\omega+1}$ implies that there is no 
$(\aleph_n,\aleph_\omega)$-Jensen matrix for any $n\geq 1$. 

\section{A characterization of the $\kappa$-Freese-Nation 
\protect\\property}\label{xchar} 
The following game over a \po\ $P$ was considered in \cite{WFN,pogame}. Let 
$\kappagame{\kappa}{P}$ be the following game 
played by Players I and II: 
in a play in $\kappagame{\kappa}{P}$, Players I and II choose 
subsets 
$X_\alpha$ and $Y_\alpha$ of $P$ of cardinality less than $\kappa$ 
alternately for 
$\alpha<\kappa$ \st
\[ 
X_0\subseteq Y_0\subseteq X_1\subseteq Y_1\subseteq\cdots\subseteq X_\alpha
\subseteq Y_\alpha\subseteq\cdots\subseteq X_\beta\subseteq 
Y_\beta\subseteq\cdots 
\]\noindent
for $\alpha\leq\beta<\kappa$. Thus a play in $\kappagame{\kappa}{P}$ 
looks like 
\[ 
\begin{array}{l@{}l@{}@{}l@{}@{}l@{}l}
\mbox{\it Player I}\ &:\ \ &X_0,\ &X_1,\ \ldots,\ &X_\alpha,\ \ldots\medskip\\ 
\mbox{\it Player II}\ &:\ \ &Y_0,\ &Y_1,\ \ldots,\ &Y_\alpha,\ \ldots
\end{array}
\]\noindent
where $\alpha<\kappa$. Player II wins the play if 
$\bigcup_{\alpha<\kappa}X_\alpha=\bigcup_{\alpha<\kappa}Y_\alpha$ is a 
$\kappa$-sub-ordering of $P$. 
Let us call a strategy $\tau$ for Player II {\em simple} if, in $\tau$, each 
$Y_\alpha$ is decided from the information of the set 
$X_\alpha\subseteq P$ alone (i.e.\ also independent of $\alpha$). 

Another notion we need here is the following generalization of 
$V_\kappa$-likeness. Let $\kappa$ be regular and $\chi$ be sufficiently 
large. For $\calD\subseteq\setof{M\prec\calH(\chi)}{\cardof{M}<\kappa}$, 
we say that $M\in[\calH(\chi)]^\kappa$ is {\em$\calD$-approachable} if 
there is an increasing  sequence $(D_\alpha)_{\alpha<\kappa}$ of elements 
of $\calD$ 
\st\medskip\\ 
\assert{a} $D_\alpha\cup\smallsetof{D_\alpha}\subseteq D_{\alpha+1}$ for 
every $\alpha<\kappa$; and \\
\assert{b} $M=\bigcup_{\alpha<\kappa}D_\alpha$. \medskip\\
Clearly $M\prec\calH(\chi)$ is $V_\kappa$-like if and only if $M$ is 
$\calD$-approachable for 
$\calD=\setof{M\prec\calH(\chi)}{\cardof{M}<\kappa}$. 

A slightly weaker version of the following theorem was announced in 
\cite{WFN}: 
\begin{Thm}
\label{correction}
Let $\kappa$ be  a regular uncountable cardinal and $\kappa\leq\lambda$. 
Suppose that 
\smallskip\\
\assert{i}
$([\mu]^{<\kappa},\subseteq)$ has a cofinal subset of cardinality $\mu$ 
for every $\mu$ \st\ $\kappa<\mu<\lambda$ and $\cf(\mu)\geq\kappa$\,; 
and\\
\assert{ii}
$\Box^{***}_{\kappa,\mu}$ holds for every $\mu$ \st\ 
$\kappa\leq\mu<\lambda$ and 
$\cf(\mu)<\kappa$.\smallskip\\
Then, for a \po\ $P$ of cardinality $\leq\lambda$, 
\tfae:\medskip\\ 
\assert{1} $P$ has the $\kappa$-FN\,;\smallskip\\
\assert{2} Player II has a simple winning strategy in 
$\kappagame{\kappa}{P}$\,;\smallskip\\
\assert{3} for some, or equivalently any sufficiently large $\chi$, and 
any $V_\kappa$-like $M\prec\calH(\chi)$ with $P$, $\kappa\in M$, we have 
$P\cap M\leq_\kappa P$\,;\smallskip\\
\assert{4} for some, or equivalently any sufficiently large $\chi$, 
there is $\calD\subseteq[\calH(\chi)]^{<\kappa}$ \st\ 
$\calD$ is cofinal in $[\calH(\chi)]^{<\kappa}$ and for any 
$\calD$-approachable $M\subseteq\calH(\chi)$, we 
have $P\cap M\leq_\kappa P$. 
\end{Thm}
Note that $\neg 0^\#$ implies the conditions \assertof{i} and 
\assertof{ii}. Also note that, for every $\lambda<\kappa^{+\omega}$, the 
condition \assertof{i} holds in ZFC. Hence the characterization above 
holds for partial orderings of cardinality $\leq\kappa^{+\omega}$ 
without any additional assumptions. \\ 
\prf
A proof of 
\assertof{1}\implies\assertof{2}\implies\assertof{3} is given in 
\cite{WFN}.
For \assertof{3}\implies\assertof{4}, suppose that $P$ 
satisfies \assertof{3}. Then $P$ together with 
$\calD=\setof{M\prec\calH(\chi)}{\cardof{M}<\kappa,\,P,\kappa\in M}$ 
satisfies \assertof{4}. 
The proof of 
\assertof{4}\implies\assertof{1} is done by induction on 
$\nu=\cardof{P}\leq\mu$. If 
$\nu\leq\kappa$, then $P$ has the $\kappa$-FN by \Lemmaof{baka}. 
For $\nu>\kappa$, let $P$ and $\calD$ be as in \assertof{4} and assume that 
\assertof{4}\implies\assertof{1} holds for every \po\ of 
cardinality $<\mu$. 
We need the following claims:
\begin{Claim}\label{claim-1}
Let $\chi^*$ be sufficiently large above $\chi$. Suppose that $M$ is an 
elementary submodel of $\calH(\chi^*)$ \st\ $P$, $\calH(\chi)$, $\calD\in M$, 
$\kappa+1\subseteq M$ and $[M]^{<\kappa}\cap M$ is cofinal in 
$[M]^{<\kappa}$ \wrt\ $\subseteq$. Then we have $P\cap M\leq_\kappa P$. 
\end{Claim}
\prfofClaim
Suppose not. then there is $b\in P$ \st\ either 
\medskip\\
\assert{a} 
$(P\cap M)\restr b$ has no cofinal subset of cardinality $<\kappa$\,; or 
\\
\assert{b} $(P\cap M)\uparrow b$ has no coinitial subset of cardinality 
$<\kappa$. 
\medskip\\
To be definite, let us assume that we have the case 
\assertof{a} --- for the case \assertof{b}, just the same argument 
will do. 
We can 
construct an increasing sequence $(N_\alpha)_{\alpha<\kappa}$ of 
elements of $\calD$ \st\medskip\\
\assert{c} $N_\alpha\in M$ and $\cardof{N_\alpha}<\kappa$ for 
$\alpha<\kappa$ (since $\kappa+1\subseteq M$, it follows that 
$N_\alpha\subseteq M$)\,;\\
\assert{d} $N_\alpha\in N_{\alpha+1}$ for every $\alpha<\kappa$\,;\\
\assert{e} $(P\cap N_\alpha)\restr b$ is not cofinal in 
$(P\cap N_{\alpha+1})\restr b$ for every $\alpha<\kappa$.\medskip\\
Then $N=\bigcup_{\alpha<\kappa}N_\alpha$ is $\calD$-approachable elementary 
submodel of $\calH(\chi)$ by \assertof{c} and \assertof{d}. 
Hence, by \assertof{4}, we have
$P\cap M\leq_\kappa P$. But, 
by \assertof{e}, 
$(P\cap N)\restr b$ has no cofinal subset of cardinality $<\kappa$. This 
is a contradiction.

To see that the construction of $N_\gamma$ is possible 
at a limit $\gamma<\kappa$, assume that 
$N_\alpha$, $\alpha<\gamma$ have been constructed in accordance with  
\assertof{c},\assertof{d} and \assertof{e}. Let 
$N'=\bigcup_{\alpha<\gamma}N_\alpha$. By \assertof{e}, we have 
$N'\subseteq M$ and $\cardof{N'}<\kappa$. Since $[M]^{<\kappa}\cap M$ is 
cofinal in $[M]^{<\kappa}$, there is some $N''\in M$ \st\ 
$N'\subseteq N''\subseteq \calH(\chi)$ and $\cardof{N''}<\kappa$. Hence 
by elementarity of $M$ and by $\calD\in M$ there is $N'''\in M\cap\calD$ \st\ 
$N''\subseteq N'''$. 
Clearly, we may let $N_\gamma=N'''$. 
For the construction at a successor step, assume that $N_\alpha$ have 
been chosen in accordance with \assertof{c},\assertof{d} and \assertof{e}. By 
assumption, there is $c\in (P\cap M)\restr b$ \st\ there is no 
$c'\in (P\cap N_\alpha)\restr b$ \st\ $c\leq c'$. By elementarity of $M$ 
there is $N^*\in M\cap\calD$ \st\ 
$N_\alpha\cup\smallsetof{N_\alpha,c}\subseteq N^*$ and 
$\cardof{N^*}<\kappa$. Then $N_{\alpha+1}=N^*$ is as desired. 
\qedofClaim
\begin{Claim}\label{claim-2}
If $Q\leq_\kappa P$, then for every 
$\calD$-approachable $M\prec\calH(\chi)$ with $Q\in M$ we have 
$Q\cap M\leq_\kappa Q$. In particular, such $Q$ also satisfies the 
condition \assertof{4}.
\end{Claim}
\prfofClaim
Let 
$M\prec\calH(\chi)$ be $\calD$-approachable with $Q\in M$. By assumption, 
we have $P\cap M\leq_\kappa P$. By 
elementarity of $M$ and since $Q\in M$, we have 
$Q\cap M\leq_\kappa P\cap M$. It follows that 
$Q\cap M\leq_\kappa P$ and hence $Q\cap M\leq_\kappa Q$. 
Now, let $\calD_0=\setof{M\in\calD}{Q\in M}$. Then it is clear that 
$Q$ satisfies the condition \assertof{4} with $\calD_0$ in place of $\calD$. 
\qedofClaim
\medskip\\
Now we are ready to prove the induction steps.\\
{\bf Case I\,:} $\nu$ is a limit cardinal or $\nu=\mu^+$ with 
$\cf(\mu)\geq\kappa$. 
Let $\nu^*=\cf\nu$. Then, by \assertof{i}, we can find an increasing sequence 
$(M_\alpha)_{\alpha<\nu^*}$ of elementary submodels of $\calH(\chi)$ \st,
for every $\alpha<\nu^*$, 
$\cardof{M_\alpha}<\nu$ and $M_\alpha$ satisfies the conditions in 
\Claimof{claim-1}; and $P\subseteq\bigcup_{\alpha<\nu^*}M_\alpha$. 
By \Claimof{claim-1}, $P\cap M_\alpha\leq_\kappa P$ for every 
$\alpha<\nu^*$. For $\alpha<\nu^*$ let
\[ P_\alpha=\left\{\,
	\begin{array}{@{}l@{}l}
		P\cap M_\alpha,\quad&\mbox{if }\alpha\mbox{ is a successor},\\[\jot]
		P\cap(\bigcup_{\beta<\alpha}M_\beta),\quad&\mbox{otherwise}.
	\end{array}\right.
\]\noindent
Then $(P_\alpha)_{\alpha<\nu^*}$ is a continuously increasing sequence of 
sub-orderings of $P$ \st\ $\cardof{P_\alpha}<\nu$ for every 
$\alpha<\nu^*$ and $P=\bigcup_{\alpha<\nu^*}P_\alpha$. We have also 
$P_\alpha\leq_\kappa P$ for every $\alpha<\nu^*$: for a successor 
$\alpha<\nu^*$ this is clear. If a limit $\alpha<\nu^*$ has cofinality 
$<\kappa$ then $P_\alpha$ can be represented as the union of an increasing 
sequence of $<\kappa$ many 
$\kappa$-sub-ordering of $P$ and hence $P_\alpha\leq_\kappa P$. If 
$\alpha<\nu^*$ is a limit with cofinality $\geq\kappa$, then 
$P_\alpha=P\cap M$ where $M=\bigcup_{\beta<\alpha}M_\beta$. Now it is 
clear that $M$ satisfies the conditions in \Claimof{claim-1}. Hence we 
again obtain that $P_\alpha=P\cap M\leq_\kappa P$. Now, by 
\Claimof{claim-2}, each of $P_\alpha$, $\alpha<\nu^*$ satisfies the 
condition \assertof{4} and hence, by induction hypothesis, has the 
$\kappa$-FN. Thus, by \Lemmaof{conti-chain}, $P$ also has the 
$\kappa$-FN. \medskip\\
{\bf Case II\,:} $\nu=\mu^+$ with $\cf(\mu)<\kappa$. Let 
$\mu^*=\cf(\mu)$. \Wolog\ we may assume 
that the underlying set of $P$ is $\nu$. By \Thmof{Jensen-exists}, there is a
$(\kappa,\mu)$-Jensen matrix 
$(M_{\alpha,\beta})_{\alpha<\nu,\beta<\mu^*}$ over $(P,\calH(\chi))$. 
For $\alpha<\nu$ and $\beta<\mu^*$, let 
$P_{\alpha,\beta}=P\cap M_{\alpha,\beta}$ and 
$P_\alpha=\bigcup_{\beta<\mu^*}P_{\alpha,\beta}$. By \assertof{j4}, the 
sequence $(P_\alpha)_{\alpha<\nu}$ is continuously increasing and 
$\bigcup_{\alpha<\nu}P_\alpha=P$. 
$\cardof{P_\alpha}\leq\mu$ for every $\alpha<\nu$ by \assertof{j1}. 
\begin{Claim}\label{claim-4} $P_\alpha\leq_\kappa P$ for every 
$\alpha<\nu$.
\end{Claim}
\prfofClaim
For $\alpha<\nu$ \st\ $\cf(\alpha)\geq\kappa$, we have 
$P_{\alpha,\beta}\leq_\kappa P$ for every sufficiently large 
$\beta<\mu^*$ by \assertof{j3} and \Claimof{claim-1}. 
Since $\mu^*<\kappa$, 
it follows that $P_\alpha\leq_\kappa P$. If $\cf(\alpha)<\kappa$, then, 
by the argument above, 
we have $P_{\alpha'}\leq_\kappa P$ for every $\alpha'<\alpha$ with 
$\cf(\alpha')\geq\kappa$. Since $P_\alpha$ can be represented as the 
union of $<\kappa$ many of such $P_{\alpha'}$'s, it follows again that 
$P_\alpha\leq_\kappa P$.
\qedofClaim
\medskip\\
Now, by \Claimof{claim-2}, each of $P_\alpha$, $\alpha<\nu$ satisfies the 
condition of \assertof{4}. Hence, by induction hypothesis, they have the 
$\kappa$-FN. By \Lemmaof{conti-chain}, it follows that $P$ also 
has the $\kappa$-FN.\qedofThm
\qedskip\\
\begin{Cor}\label{appl-of-char}
Suppose that $\kappa$ and $\lambda$ satisfy \assertof{i}, \assertof{ii} 
in \Thmof{correction} and $2^{<\kappa}=\kappa$. Then:\smallskip\\
\assert{a} Every $\kappa$-cc \cBa\ of cardinality $\leq\lambda$ has the 
$\kappa$-FN. \smallskip\\
\assert{b} For any $\mu$ \st\ $\mu^{<\kappa}\leq\lambda$, the \po\ 
$([\mu]^{<\kappa},\subseteq)$ has the $\kappa$-FN.
\end{Cor}
\prf
 Let $\chi$ be sufficiently large. 
For \assertof{a}, 
let $B$ be a $\kappa$-cc \cBa. We show that $B$ satisfies \assertof{3} in 
\Thmof{correction}. Let 
$M\prec\calH(\chi)$ be $V_\kappa$-like with $B$, $\kappa\in M$. 
By $2^{<\kappa}=\kappa$, \Lemmaof{easy} and by the remark after the lemma, 
we have $[M]^{<\kappa}\subseteq M$. Hence 
$B\cap M$ is a complete subalgebra of $B$. It follows that 
$B\cap M\leq_\kappa M$. 
For \assertof{b}, let 
$M\prec\calH(\chi)$ be $V_\kappa$-like with $\lambda$, $\kappa\in M$. 
Then as above we have $[M]^{<\kappa}\subseteq M$. Hence, 
letting $X=\lambda\cap M$, we have 
$[\lambda]^{<\kappa}\cap M=[X]^{<\kappa}$. It follows that 
$[\lambda]^{<\kappa}\cap M\leq_\kappa [\lambda]^{<\kappa}$. 
\qedofCor
\section{Chang's Conjecture for $\aleph_\omega$}\label{chang}
Recall that 
$(\kappa,\lambda)\pfeil(\mu,\nu)$ is the following assertion: 
\begin{assertion}{}\it
For any structure ${\cal A}=(A,U,\ldots)$ of countable signature with 
$\cardof{A}=\kappa$, $U\subseteq A$ and $\cardof{U}=\nu$, there is an 
elementary substructure ${\cal A}'=(A',U',\ldots)$ of $\cal A$ \st\ 
$\cardof{A'}=\mu$ and $\cardof{U'}=\nu$. 
\end{assertion}
In \cite{levinski}, a model of ZFC $+$ GCH $+$ Chang's Conjecture for 
$\aleph_\omega$, i.e.\ 
$(\aleph_{\omega+1},\aleph_\omega)\pfeil(\aleph_1,\aleph_0)$, is 
constructed starting from a model with a cardinal having a property 
slightly stronger than huge. The following theorem together with 
\Corof{appl-of-char} shows that the $\aleph_1$-FN of the \po\ 
$([\aleph_\omega]^{\aleph_0},\subseteq)$ is independent from ZFC 
(or even from ZFC $+$ GCH). 
\begin{Thm}\label{ccchang}
Suppose that $(\aleph_\omega)^{\aleph_0}=\aleph_{\omega+1}$ and 
$(\aleph_{\omega+1},\aleph_\omega)\pfeil(\aleph_1,\aleph_0)$. Then 
$([\aleph_\omega]^{\aleph_0},\subseteq)$ does not have the $\aleph_1$-FN.
\end{Thm}
\prf
\nc{\alo}{{\aleph_{\omega}}}
\nc{\aloo}{{\aleph_{\omega+1}}}
Assume to the contrary that there is an $\aleph_1$-FN mapping 
$\mapping{F}{[\alo]^{\aleph_0}}%
	{[[\alo]^{\aleph_0}]^{\aleph_0}}$. 
Let us fix an enumeration $(b_\alpha)_{\alpha<\aloo}$ of 
$[\alo]^{\aleph_0}$ and consider the structure:
\[ \calA=(\aloo, \alo,\leq, E,f,g,h),
\]\noindent
where\medskip\\
\assert{1} $\leq$ is the canonical ordering on $\aloo$;\smallskip\\
\assert{2} 
$E=\setof{(\alpha,\beta)}{\alpha\in\alo,\beta\in\aloo,\alpha\in b_\beta}$;
\smallskip\\
\assert{3} $\mapping{f}{\aloo\times\aloo}{\aloo}$ is \st, for each 
$\alpha\in\aloo$,
\[ F(b_\alpha)=\setof{b_{f(\alpha,n)}}{n\in\omega};
\]\noindent
\assert{4} $\mapping{g}{\aloo\times\aloo}{\alo}$ is \st, for each 
$\alpha\in\aloo$, $g(\alpha,\cdot)\restr\alpha$ is an injective mapping 
from $\alpha$ to $\alo$; and
\smallskip\\
\assert{5} $\mapping{h}{\aloo\times\aloo\times\aloo}{\omega+1}$ is \st\ 
for each $\alpha$, $\beta\in\aleph_{\omega+1}$, 
$h(\alpha,\beta,\cdot)\restr(b_\alpha\cap b_\beta)$ is injective.\medskip\\
Note that, by \assertof{5} and since $\omega$ is definable in $\calA$, we 
can express ``$b_\alpha\cap b_\beta$ is 
finite'' as a formula $\varphi(\alpha,\beta)$ in the language of $\calA$. 
Now, by assumption there is an 
elementary substructure $\calA'=(A',U',\leq',E',f',g',h')$ of $\calA$ 
\st\ $\cardof{A'}=\aleph_1$ and $\cardof{U'}=\aleph_0$. 
\begin{Claim}\label{claim=1}
$otp(A')=\omega_1$.
\end{Claim}
\prfofClaim
By \assertof{4} and elementarity of $\calA'$, every initial segment of 
$A'$ can be mapped into $U'$ injectively and hence countable. Since 
$\cardof{A'}=\aleph_1$, it follows that $otp(A')=\omega_1$.
\qedofClaim
\begin{Claim}\label{claim=2}
For any $\alpha<\aloo$, there is $\gamma<\aloo$ \st\ 
$b_\beta\cap b_\gamma$ is finite for every $\beta<\alpha$.
\end{Claim}
\prfofClaim
Since $\cardof{\alpha}\leq\alo$, we can find a partition 
$(I_n)_{n\in\omega}$ of $\alpha$ \st\ $\cardof{I_n}<\alo$ for every 
$n<\omega$. For $n<\omega$, let 
$\eta_n=\min(\alo\setminus\bigcup\setof{b_\xi}{\xi\in\bigcup_{m\leq n}I_m})$. 
Let $z=\setof{\eta_n}{n\in\omega}$ and $\gamma<\aloo$ be \st\ 
$b_\gamma=z$. For any $\beta<\alpha$, if $\beta\in I_{m_0}$ for some 
$m_0<\omega$, then we have 
$b_\beta\cap b_\gamma\subseteq\setof{\eta_n}{n<m_0}$. Thus this $\gamma$ 
is as desired.
\qedofClaim 
\begin{Claim}\label{claim=3}
For any countable $I\subseteq A'$, there is $\gamma\in A'$ \st\ 
$b_\beta\cap b_\gamma$ is finite for every $\beta\in I$.
\end{Claim}
\prfofClaim
By \Claimof{claim=1}, there is $\alpha\in A'$ \st\ $I\subseteq\alpha$. By 
elementarity of $\calA'$, the formula with the parameter $\alpha$ 
expressing the assertion of \Claimof{claim=2} for this $\alpha$ holds in 
$\calA'$. Hence there is some $\gamma\in A'$ \st\ 
$b_\beta\cap b_\gamma$ is finite for every $\beta\in A'\cap\alpha$. 
\qedofClaim 
\medskip\\
Let 
\[ I=\setof{\xi\in A'}{b_\xi\in F(U')}. 
\]\noindent
Then $I$ is countable. Hence, by 
\Claimof{claim=3}, there is $\gamma\in A'$ \st\ $b_\beta\cap b_\gamma$ is 
finite for every $\beta\in I$. As $b_\gamma\subseteq U'$ (this holds in 
virtue of $h(\gamma,\gamma,\cdot)$), there is 
$b\in F(b_\gamma)\cap F(U')$ \st\ $b_\gamma\subseteq b\subseteq U'$. Let 
$b=b_{\xi_0}$. Then $\xi_0\in I$ and 
$\cardof{b_\gamma\cap b_{\xi_0}}=\cardof{b_\gamma}=\aleph_0$. This is a 
contradiction to the choice of $\gamma$. 
\qedofThm
\qedskip
\par
We do not know if the assumption of \Thmof{ccchang} yields a counter 
example to \Corof{appl-of-char},\,\assertof{a} for $\kappa=\aleph_1$. Or, 
more generally:
\begin{Problem}
Is there a model of $\ZFC+\GCH$ where some ccc \cBa\ does not have the 
$\aleph_1$-FN\,? 
\end{Problem}
Of course, we need the consistency strength of some large cardinal to 
obtain such a model by \Corof{appl-of-char}. 
\par The following corollary slightly improves Theorem 4.1 in 
\cite{foreman-magidor}. 
\begin{Cor}
Suppose that $(\aleph_\omega)^{\aleph_0}=\aleph_{\omega+1}$ and 
$(\aleph_{\omega+1},\aleph_\omega)\pfeil(\aleph_1,\aleph_0)$. Then 
the equivalence of the assertions in \Thmof{correction} fails. Hence we 
have $\neg \Box^{***}_{\aleph_1,\aleph_\omega}$ under these assumptions. 
\end{Cor}
\prf 
By \Thmof{ccchang} and \Corof{appl-of-char},\,\assertof{b}. 
\qedofCor\qedskip\\
Similarly we can prove $\neg\Box^{***}_{\aleph_n\aleph_\omega}$ for every 
$n\in\omega$ under the assumption of 
$2^{\aleph_\omega}=\aleph_{\omega+1}$ and 
$(\aleph_{\omega+1},\aleph_\omega)\pfeil(\aleph_1,\aleph_0)$. 

\section{Cohen models}\label{cohen}
\newcommand{\dotleq}{\mathrel{\dot{\leq}}}
\newcommand{\dotleqd}{\mathrel{\dot{\leq}'}}
Let $V$ be our ground model and let 
$G$ be a $V$-generic filter over $P=\Fn(\tau,2)$ for some $\tau$. Suppose 
that $B$ is a ccc 
\cBa\ in $V[G]$. \Wolog\ we may assume that the underlying set 
of $B$ is a set $X$ in $V$. $B$ is said to be {\it tame} if there is a 
$P$-name $\dotleq$ of partial ordering of $B$ and a mapping 
$\mapping{t}{X}{[\tau]^{\aleph_0}}$ in $V$ \st, for every $p\in P$ and $x$, 
$y\in X$, if $p\forces{P}{x\dotleq y}$, then 
$p\restr(t(x)\cup t(y))\forces{P}{x\dotleq y}$. 
A lot of `natural' ccc \cBas\ in $V[G]$ are contained in the class of tame 
\Bas: 
\begin{Lemma}
Let $G$ be as above. Suppose that 
$V[G]\models{B\xmbox{ is a ccc complete \Ba}}$ 
and either:\smallskip\\
\assert{i} there is a \Ba\ $B'$ in $V$ \st\ $B'$ is dense subalgebra of 
$B$ in $V[G]$; or\\ 
\assert{ii} $B$ is the completion of a Suslin forcing in $V[G]$.
\smallskip\\
Then $B$ is tame.\qed
\end{Lemma}
For Suslin forcing, see e.g.\ \cite{bartoszynski-book}.
\begin{Thm}\label{cohen-model-rev}
Let $P$, $G$ be as above and $\lambda$ an infinite cardinal. 
Assume that, in $V$, \smallskip\\
\assert{i} $\mu^{\aleph_0}=\mu$ for every regular uncountable $\mu$ \st\ 
$\mu<\lambda$; and\\
\assert{ii} $\Box^{***}_{\aleph_1,\mu}$ holds for every $\mu$ \st\ $\kappa\leq\mu<\lambda$ and 
$\cf(\mu)=\omega$.\smallskip\\
Then, for any tame ccc \cBa\ $B$ in $V[G]$ of cardinality 
$\leq\lambda$ we have $V[G]\models{B\mbox{ has the }\aleph_1\mbox{-FN}}$. 
\end{Thm} 
\prf Let $X$, $\dot\leq$ and $t$ be as in the definition of tameness for 
$B$ above. 
We may assume 
$\forces{P}{X\xmbox{ is the underlying set of }\dot{B}}$ and  
$\forces{P}{\dot{B}\xmbox{ is a ccc \cBa}}$ where $\dot{B}$ is a 
$P$-name for $B$. 
Let $\chi$ be sufficiently large. The following is the key lemma to 
the proof:
\begin{Claim}\label{claim-x}
Suppose that $M\prec\calH(\chi)$ is \st\ $\tau$, $X$, $\dotleq$, $t\in M$ 
and $[M]^{\aleph_0}\subseteq M$. Let $I=\tau\cap M$, $P'=\Fn(I,2)$, 
$G'=G\cap P'$,
$X'=X\cap M$ and $\dotleqd={\dotleq}\cap M$. Then\smallskip\\
\assert{a} $\dotleqd$ is a $P'$-name,
$\forces{P}{B'=(X',{\dotleqd})\mbox{ is a ccc \cBa}}$,
$B'$ is tame (in $V[G']$) 
and the infinite sum $\Sigma^B$ in $V[G]$ is an extension of the infinite 
sum $\Sigma^{B'}$ in $V[G']$.\smallskip\\
\assert{b} $\forces{P}{(X',{\dotleqd})\leq_{\sigma}(X,{\dotleq})}$.
\end{Claim}
\prfofClaim
\assertof{a}: It is easy to see that 
$\forces{P'}{(X',\dotleqd)\xmbox{ is a \Ba}}$ and 
$\forces{P}{{(X',\dotleqd)}\xmbox{ is a subalgebra of }(X,\dotleq)}$. Since 
$\forces{P}{(X,\dotleq)\xmbox{ has the ccc}}$, it follows that 
$\forces{P'}{{(X',\dotleqd)}\xmbox{ has the ccc}}$. 
By elementarity of $M$, it is also easy to see that $t\restr X'$ witnesses 
the tameness of $B'$ in $V[G']$. 
To see that 
$\forces{P'}{(X',\dotleqd)\xmbox{ is complete}}$ it is enough to see that 
every countable subset of $B'$ has its supremum in $V[G']$. Let $\dot{C}$ 
be a $P'$-name of countable subset of $X'$. \Wolog, we may assume that 
$\dot{C}$ is countable and consists of sets of the form $(p,\check{x})$ 
where $p\in P'$ and $x\in X'$. Since $[M]^{\aleph_0}\subseteq M$, 
$\dot{C}\in M$. Clearly, we have
$M\models{\dot{C}\xmbox{ is a }P\xmbox{-name of countable subset of }X}$.
Hence, 
$M\models{\exists p\in P\ \exists y\in X(p\forces{P}{\Sigma\dot{C}=y})}$. 
Let $p\in P$ and $y\in X$ be such elements of $M$. Then $p\in P'$ and 
$y\in X'$. By elementarity of $M$, we have 
$p\forces{P}{\Sigma^B\dot{C}=y}$. On the other hand, from 
$M\models{p\forces{P}{\Sigma\dot{C}=y}}$ it follows that 
$p\forces{P'}{\Sigma^{B'}\dot{C}=y}$. 
\medskip\\
\assertof{b}: 
By assumption, for $x\in X$ and $y\in X'$, we have $y\leq x$ 
in $V[G]$, if and only if there is $p\in G$, 
$\dom(p)\subseteq t(x)\cup t(y)$ \st\ $p\forces{P}{y\dotleq x}$. For 
$q\in G\cap \Fn(t(y),2)$, the set 
$U_q=\setof{y\in X'}{\exists p\in G'\,(p\cup q\forces{P}{y\dotleq x})}$ 
is in $V[G']$. Hence, by \assertof{a}, $\Sigma^{B}U_q$ is an element in 
$B'$. Since $X'$ is closed under $t$, it follows that 
$\setof{\Sigma^{B}U_q}{q\in G\cap \Fn(t(y),2)}$ is cofinal in 
$B'\restr y$. 
\qedofClaim\medskip\par
Now, let $\nu=\cardof{X}$. \Wolog, we may assume that $X=\nu$. We show by 
induction on $\nu$ that 
$\forces{P}{\dotB\xmbox{ has the }\aleph_1\xmbox{-FN}}$. For 
$\nu\leq\aleph_1$ the assertion follows from \Lemmaof{baka} (applied in 
$V^P$). In the induction steps, we mimic the proof of \Thmof{correction}. 
Let $\chi$ be sufficiently large.\medskip\\ 
{\bf Case I\,:} 
$\nu$ is a limit cardinal or $\nu=\mu^+$ for some $\mu$ with 
$\cf(\mu)\geq\omega_1$. By \assertof{i}, we can construct a continuously 
increasing sequence $(M_\alpha)_{\alpha<\nu}$ of elementary submodels of 
$\calH(\chi)$ \st\
\begin{assertion}{}\mbox{}%
\lassert{0} $\tau$, $\nu$, $\dotleq$, $t\in M_0$;\\
\lassert{1} $\cardof{M_\alpha}<\nu$ for every $\alpha<\nu$;\\
\lassert{2} $[M_{\alpha+1}]^{\aleph_0}\subseteq M_{\alpha+1}$ for every 
$\alpha<\nu$ (note that it follows that the inclusion also holds for 
every limit $<\nu$ of cofinality $\geq\omega_1$); and\\
\lassert{3} $\nu\subseteq\bigcup_{\alpha<\nu} M_\alpha$. 
\end{assertion}
For each $\alpha<\nu$, let $X_\alpha=X\cap M_\alpha$ and 
$\dotleq^\alpha=\dotleq\cap M_\alpha$ and let $\dotB_\alpha$ be the $P$-name 
corresponding to $(X_\alpha, \dotleq^\alpha)$. By \Claimof{claim-x}, we 
have 
$\forces{P}{\dotB_\alpha\xmbox{ is ccc \cBa\ and }
	\dotB_\alpha\leq_{\aleph_1}\dotB}$, for all $\alpha<\nu$ \st\ 
either $\alpha$ is a successor or of cofinality $\geq\omega_1$. By 
induction hypothesis, we have 
$\forces{P}{\dotB_\alpha\xmbox{ has the }\aleph_1\xmbox{-FN}}$ for such 
$\alpha$'s. Hence by \Lemmaof{conti-chain} and the remark after the lemma 
(applied in $V^P$) it follows that 
$\forces{P}{\dotB\xmbox{ has the }\aleph_1\xmbox{-FN}}$. \medskip\\
{\bf Case II\,:} $\nu=\mu^+$ for a $\mu$ with $\cf(\mu)=\omega$. By 
\assertof{ii}, there is an $(\aleph_1,\mu)$-Jensen matrix 
$(M_{\alpha,n})_{\alpha<\nu,n<\omega}$ over $(\tau,\nu,\dotleq,t)$. For 
$\alpha<\nu$, let $M_\alpha=\bigcup_{n<\omega}M_{\alpha,n}$. 
For $\alpha<\nu$ and $n<\omega$, let 
$X_{\alpha,n}=X\cap M_{\alpha,n}$, 
$\dotleq^{\alpha,n}=\dotleq\cap M_{\alpha,n}$ and $\dotB_{\alpha,n}$ be 
the $P$-name corresponding to $(X_{\alpha,n},\dotleq^{\alpha,n})$. 
Likewise, let $X_{\alpha}=X\cap M_{\alpha}$, 
$\dotleq^{\alpha}=\dotleq\cap M_{\alpha}$ and $\dotB_{\alpha}$ be 
the $P$-name corresponding to $(X_{\alpha},\dotleq^{\alpha})$. 
Then we have $X_\alpha=\bigcup_{n<\omega}X_{\alpha,n}$, 
$\dotleq^\alpha=\bigcup_{n<\omega}\dotleq^{\alpha,n}$ and 
$\forces{P}{\dotB_\alpha=\bigcup_{n<\omega}\dotB_{\alpha,n}}$. By 
\Lemmaof{easy} and \assertof{i}, \assertof{j3'} holds for the Jensen 
matrix. Hence, by \Claimof{claim-x}, we have 
$\forces{P}{\dotB_{\alpha,n}\leq_{\aleph_1}\dotB_\alpha\xmbox{ and }
	\dotB_{\alpha,n}\xmbox{ is a ccc \cBa}}$ for every $\alpha<\nu$ with 
$\cf(\alpha)>\omega$ and every sufficiently large $n<\omega$. By 
induction hypothesis, it follows that 
$\forces{P}{\dotB_{\alpha,n}\xmbox{ has the }\aleph_1\xmbox{-FN}}$ for 
such $\alpha$ and $n$. By \Lemmaof{conti-chain} (applied in $V^P$) it 
follows that 
$\forces{P}{\dotB_\alpha\xmbox{ has the }\aleph_1\xmbox{-FN}}$ for every 
$\alpha<\nu$ with $\cf(\alpha)>\omega$. Hence again by 
\Lemmaof{conti-chain} and the remark after that (applied in $V^P$) we 
obtain that $\forces{P}{\dotB\xmbox{ has the }\aleph_1\xmbox{-FN}}$. 
\qedofThm
\begin{Cor}
Suppose that $V=L$ holds and let $G$ be $V$-generic over $P=\Fn(\tau,2)$ 
for some $\tau$. Then (among others) the following ccc \cBas\ have the 
$\aleph_1$-FN: $\calC_\kappa$ ($\cong$ the completion of $\Fn(\kappa,2)$) for 
any $\kappa$; $\psof{\omega}$ (hence also $\psof{\omega}/{\rm fin}$); 
$Borel\/(\reals)/\mbox{Null-sets}$.\qed
\end{Cor}
In connection with \Thmof{cohen-model-rev}, we would like to mention the 
following open problems:
\begin{Problem}
Assume that $V[G]$ is a model obtained by adding Cohen reals to a model 
of $\ZFC +\CH$. Is it true that $\powersetof{\omega}$ has the 
$\aleph_1$-FN in $V[G]$\,?
\end{Problem}
By \Thmof{cohen-model-rev}, the answer to this question is positive if 
the number of added Cohen reals is less than $\aleph_\omega$. 
\begin{Problem}
Does \Thmof{cohen-model-rev} hold also without the assumption of 
tameness\,? 
\end{Problem}
Or, more generally:
\begin{Problem}
Are the following equivalent?\\
\assert{i} $\powersetof{\omega}$ has the $\aleph_1$-FN;\\
\assert{ii} every ccc \cBa\ has the $\aleph_1$-FN.
\end{Problem}
\section{Lusin gap}\label{lusin}
For a regular $\kappa$, an almost disjoint family 
$\calA\subseteq[\omega]^{\aleph_0}$ is said to be a {\em $\kappa$-Lusin 
gap} if 
$\cardof{\calA}=\kappa$ and for any $x\in[\omega]^{\aleph_0}$ either 
$\cardof{\setof{a\in\calA}{\cardof{a\setminus x}<\aleph_0}}<\kappa$ or 
$\cardof{\setof{a\in\calA}{\cardof{a\cap x}<\aleph_0}}<\kappa$. 
\begin{Thm}\label{lusin-gap}
Assume that $\psof{\omega}$ has the $\aleph_1$-FN. Then there is no 
$\aleph_2$-Lusin gap.
\end{Thm}
\prf
Let $\mapping{f}{\psof{\omega}}{[\psof{\omega}]^{\aleph_0}}$ be an 
$\aleph_1$-FN mapping. We may assume that 
$f(a)=f(b)=f(\omega\setminus b)$ for all $a$, $b\in\psof{\omega}$ with 
$a=^*b$. Thus $x\subseteq^*y$ implies that there is 
$z\in f(x)\cap f(y)$ \st\ $x\subseteq^* z\subseteq^* y$ and 
$\cardof{x\cap y}<\aleph_0$ implies that there is $z\in f(x)\cap f(y)$ 
\st\ $x\subseteq^* z$ and $\cardof{z\cap y}<\aleph_0$. 

Suppose that $\calA=\setof{a_\alpha}{\alpha<\omega_2}$ is an almost 
disjoint family of subsets of $\omega$. We show that $\calA$ is not an 
$\aleph_2$-Lusin gap. 
Let $\chi$ be sufficiently large 
regular cardinal and consider the model 
$\calH=(\calH(\chi),\in,\unlhd)$ where $\unlhd$ is any 
well-ordering on $\calH$. Let $N$ be an elementary submodel of $\calH$ 
\st\ $\calA$, $f\in N$, $N\cap\omega_2\in\omega_2$ and 
$\cf(\delta)=\omega_1$ for $\delta=N\cap\omega_2$. For $\alpha\in N$ we 
have $\cardof{a_\alpha\cap a_\delta}<\aleph_0$ and hence 
$a_\alpha\subseteq^*(\omega\setminus a_\delta)$. Thus there is 
$b_\alpha\in f(a_\alpha)\cap f(a_\delta)$ \st\ 
$a_\alpha\subseteq^*b_\alpha\subseteq^*(\omega\setminus a_\delta)$. Since 
$f(a_\delta)$ is countable and $\cf(\delta)=\omega_1$ there is 
$b\in f(a_\delta)$ \st\ $I=\setof{\alpha<\delta}{b_\alpha=b}$ is cofinal 
in $\delta$. We show that $b$ witnesses that $\calA$ is not an 
$\aleph_2$-Lusin gap, i.e., 
$J=\setof{\alpha<\omega_2}{a_\alpha\subseteq^* b}$ and 
$K=\setof{\alpha<\omega_2}{\cardof{b\cap a_\alpha}<\aleph_0}$ both have 
cardinality $\aleph_2$. 
\begin{Claim}
$\cardof{J}=\aleph_2$.
\end{Claim}
\prfofClaim
First note that $b\in N$ since $b\in f(a_\alpha)$ for any 
$\alpha\in I\subseteq N$. Hence we have $J\in N$ and $I\subseteq J$. 
Since $I$ 
is cofinal in $N\cap\omega_2$, we have 
$N\models{J\mbox{ is cofinal in }\omega_2}$. By 
elementarity it follows that 
$\calH\models{J\mbox{ is cofinal in }\omega_2}$. Hence 
$J$ is really cofinal in $\omega_2$.
\qedofClaim
\begin{Claim}
$\cardof{K}=\aleph_2$.
\end{Claim}
\prfofClaim
Since $b\in N$ it follows that $K\in N$. For 
$\beta\in N\cap\omega_2=\delta$, we have 
$\calH\models{\delta\in K\wedge\beta<\delta}$. Hence 
$\calH\models{K\not\subseteq\beta}$ and 
$N\models{K\not\subseteq\beta}$ by elementarity. It follows that 
$N\models{K\xmbox{ is unbounded in }\omega_2}$. Hence, again by 
elementarity, 
$\calH\models{K\xmbox{ is unbounded in }\omega_2}$. Thus $K$ is really 
unbounded in $\omega_2$. 
\qedofClaim\\
\qedofThm
\begin{Cor}
${\bf b}=\aleph_1$ or even the statement ``\,$\psof{\omega}$ does not 
contain any strictly increasing $\subseteq^*$-chain of length 
$\omega_2$'' does not 
imply that $\psof{\omega}$ has the $\aleph_1$-FN. 
\end{Cor}
\prf Suppose that our ground model $V$ satisfies the CH. 
Koppelberg and Shelah \cite{koppelberg-shelah} proved that the 
forcing with $\Fn(\omega_2,2)$ can be 
represented as a two step iteration $A*\dot{B}$ where 
$\forces{A}{\psof{\omega}\xmbox{ contains an }\aleph_2\xmbox{-Lusin gap}}$. 
Thus, by \Thmof{lusin-gap}, we have 
$\forces{A}{\psof{\omega}\xmbox{ does not habe the }\aleph_1\xmbox{-WF}}$. 
On the other hand, we have 
$\forces{A*\dot{B}}{\psof{\omega}\xmbox{ does have the }\aleph_1\xmbox{-WF}}$ 
by \Thmof{cohen-model-rev}. 
Hence 
$\forces{A*\dot{B}}{%
	\xmbox{there is no strictly increasing }\subseteq^*
			\xmbox{-chain in }\psof{\omega}\xmbox{ of length }\omega_2}$. 
It follows that 
$\forces{A}{%
	\xmbox{there is no strictly increasing }\subseteq^*
			\xmbox{-chain in }\psof{\omega}\xmbox{ of length }\omega_2}$. 
\qedofCor 
\newpage

\end{document}